\documentclass[11pt]{amsart}
\usepackage{amsfonts,amsmath,amsthm}
\usepackage{amssymb}
\usepackage{xcolor}
\usepackage{mathrsfs}
\usepackage{bbold}

\newcommand{\overbar}[1]{\mkern 1.5mu\overline{\mkern-1.5mu#1\mkern-1.5mu}\mkern 1.5mu}

\usepackage{fullpage}
\usepackage{xspace}

\numberwithin{equation}{section}

\theoremstyle{plain}
\newtheorem{thm}{Theorem}[section]
\newtheorem{proposition}{Proposition}[section]
\newtheorem{lemma}{Lemma}[section]
\newtheorem*{theorem*}{Theorem}

\theoremstyle{definition}
\newtheorem{defn}{Definition}[section]
\newtheorem{rem}{Remark}[section]

\newcommand{\zb}{\overline{\zeta}}

\newcommand{\wb}{\overline{w}}

\newcommand{\ttr}{\operatorname{tr}}

\newcommand{\pp}{\tilde{p}}
\newcommand{\qq}{\tilde{q}}
\newcommand{\lc}{\overline{\lambda}}

\newcommand{\cb}{\overline{c}}
\newcommand{\ab}{\overline{a}}

\title{On the M\"{o}bius invariant Principal functions of Pincus}

\author{Sagar Ghosh}
\address{Indian Statistical Institute, Bangalore}
\email{sagarghosh1729@gmail.com}

\author{Gadadhar Misra}

\address{Indian Statistical Institute, Bangalore and Indian Institute of Technology, Gandhinagar}

\email{gm@isibang.ac.in}

\thanks{The second author gratefully acknowledges the financial support from the Science and Engineering Research Board (SERB) in the form of a ~J~C~Bose National Fellowship.}


\dedicatory{This paper is dedicated to Professor Jan
Stochel on the occasion of his seventieth birth anniversary.}

\begin{document}
\begin{abstract}
    In this semi-expository short note, we prove that the only homogeneous \textit{pure} hyponormal operator $T$  with $\operatorname{rank} (T^*T-TT^*) =1$, modulo unitary equivalence, is the unilateral shift. 
\end{abstract}

\maketitle 
\section{Introduction}
In this paper, a Hilbert space $\mathcal H$ is assumed to be complex and separable and an operator $T$ on $\mathcal H$ is assumed to be linear and bounded. The algebra of bounded linear operators on a complex separable Hilbert  space $\mathcal H$ is denoted by $\mathcal L(\mathcal H)$. An operator $A\in \mathcal L(\mathcal H)$ is said to be \emph{hyponormal} if $[A^*,A]:=A^*A - AA^*$ is \emph{non-negative}, that is, $\langle [A^*,A]f, f \rangle \geqslant 0$ for all $f\in \mathcal H$.

Let $\mathcal H$ be a Hilbert space and $\{e_n\}_{n\geqslant 0}$ be an orthonormal basis in $\mathcal H$. For any bounded non-negative  operator $B$ acting on $\mathcal H$, define its trace by setting
$$
\operatorname{tr}(B)=\sum_n\langle B e_n, e_n\rangle.
$$
This definition of $\operatorname{tr}(B)$ does not depend on the choice of the orthonormal basis that was chosen to define it.

An operator $A\in \mathcal L(\mathcal H)$ is said to be in the \emph{trace class} $\mathcal S_1(\mathcal H)$ if $\sum_{n=0}^{\infty}\langle |A|e_n,e_n \rangle$ is finite. As usual, here $|A|$ is the unique positive square root of the self-adjoint operator $A^*A$. 

The $s$-numbers $\{s_j(T)\}_{j=1}^\infty$ of a compact operator $T$ are
the eigenvalues of $(T^* T)^{\tfrac{1}{2}}$, counted with multiplicity and arranged in decreasing order.  The trace norm is also given by the formula:
\[\|T\|_1 = \sum_{j=1}^\infty s_j(T).\]
Let $T$ be a trace class operator. Set 
\[\Lambda_T := \{\lambda_j(T) : j = 1,2 , \ldots \nu(T) \}\]
be an enumeration of the non-zero eigenvalues of $T$ counting 
multiplicities. The determinant of the operator $I + T$ is defined as follows: 
\[\det (I+T) = \begin{cases} 
\displaystyle\prod_{j=1}^{\nu(T)}  (1+ \lambda_j(T)), & \Lambda_T \not= \emptyset \\
\phantom{!!}1& \Lambda_T = \emptyset
\end{cases}\] 
In case $\nu(T)$ is infinite, the convergence of the product defining the determinant follows from the inequality $\sum_{j=1}^\infty \lambda_j(T) \leq \|T\|_1$, see \cite[Chapter II]{GK}. We need the following crucial relationship between the trace and the determinant.

Recall the Jacobi formula for matrix exponential, namely, $\det \exp(R) = \exp(\operatorname{tr}(B)$. Now, suppose that $T$ is a trace class operator with $\|T\|_1 < 1$. Then we define  
$\det(I+T) = \exp \operatorname{tr}(\log(I+T))$. Here $\log(I+T)$ is the logarithm of $I+T$ given by the series (convergent in the norm $\|\cdot\|_1$)
\begin{equation}\label{detdef}
\log(I+T) = -\sum_{n=1}^\infty (-1)^n \frac{T^n}{n},\end{equation}
see \cite[pp. 81]{KC}.

\begin{defn}
A natural number $m$ is said to be the (rational) multiplicity of an operator $T\in \mathcal L(\mathcal H)$ if there exist vectors $\{x_i\}_{i\in I}$, for some indexing set $I$ with $|I| = m$, such that 
\[H = \bigvee \big \{f(T)x_i, i\in I, f\in \text{Rat}(\sigma(T))\big \},\]
where $\text{Rat}(\sigma(T))$ is the set of all rational functions $r$ of the form $\frac{p}{q}$ for a pair of polynomials $p$ and $q$ with $q$ not zero on $\sigma(T)$. 
\end{defn}
In this short note we study the class of hyponormal operators $T$ with finite multiplicity. The remarkable inequality \cite{BS} of Berger and Shaw  
\begin{equation}\label{BSinq}\operatorname{tr}\left[\mathrm{T}^*, \mathrm{~T}\right] \leqslant \frac{m}{\pi} \operatorname{Area} (\sigma(\mathrm{T}))\end{equation}
ensures that  the self-commutator $[T^*,T]$ of such an operator is in  the trace class. An immediate corollary is an inequality due to Putnam \cite{CRP}: If $T\in \mathcal L(\mathcal H)$ is hyponormal, then
\begin{equation} \label{Putnam}
    \|[T^*,T]\| \leq \frac{1}{\pi}\mu(\text{Area}(T)). 
\end{equation}
The verification below of Putnam's inequality  is taken from \cite[Chapter VI, Theorem 2.1]{MP}.  
Pick a non-zero vector $x\in \mathcal H$ and set  
\begin{align*}
    \mathcal H_x: = \bigvee\big \{f(T)x ; f\in \text{Rat}(\sigma(T))\big \}.
\end{align*} Let $T_x:\mathcal H_x \to \mathcal H_x$ be the restriction of the operator $T$ to $\mathcal H_x$. The operator $T_x$ is evidently hyponormal and it is rationally cyclic of multiplicity $1$. We have 
\begin{align*}
    \langle [T^*,T]x,x]\rangle  
    & = \|Tx\|^2 - \|T^*x\|^2 \\
    & \leqslant \|T_x x\|^2 - \|T^*_x x\|^2 \\
    & = \langle [T_x^*,T_x]x,x \rangle \\
    & \leqslant \operatorname{tr}([T_x^*,T_x])\\
    & \leqslant \frac{1}{\pi}\text{Area}(\sigma(T_x))\\
    & \leqslant \frac{1}{\pi}\text{Area}(\sigma(T)),
\end{align*} where the penultimate inequality follows from Berger-Shaw inequality \eqref{BSinq} and the last inequality is a consequence of the spectral inclusion $\sigma(T_x)\subseteq \sigma(T)$. 
\begin{rem}
Among many consequences of Putnam's inequality, we single out one that we will need in what follows, namely, if $T$ is a \textit{pure} hyponormal operator, then $\text{Area}(\sigma(T)) > 0$.
\end{rem}
Moreover, we note that the \textit{determinantal formula} due to Carey and Pincus, discussed below, connects the \textit{Principal function} $g_T$ of the operator $T$ with the trace of $[T^*,T]$ using the Helton-Howe \textit{trace formula}. For a recent account, one may consult the book \cite{BP}.

\begin{defn}
The bi-holomorphic automorphism group M\"{o}b of the unit disc consists of rational functions $\varphi$  of the form: \[\varphi(z) = \beta \frac{z-a}{1-\bar{a}z},\,\, \beta\in \mathbb T \text{ and }a\in \mathbb D,\] where, $\mathbb T$ and $\mathbb D$ denote the unit circle and the open unit disc respectively.   
\end{defn}

For an opeartor $T$ with the spectrum $\sigma(T)$  contained in the closed unit disc $\overbar{\mathbb D}$, by the spectral mapping theorem, $0\not \in \sigma(I - \bar{a} T)$ for any $a\in \mathbb D$.  Hence, the operator $I - \bar{a} T$, $a\in \mathbb D$, is invertible. 
\begin{defn}
An operator $T$ with $\sigma(T) \subseteq \overbar{\mathbb D}$ is said to be \textit{homogeneous} if the operator \[\varphi(T):= \beta (T-a)(I - \bar{a} T)^{-1},\,\, \beta \in \mathbb T,\, a\in\mathbb D,\] is unitarily equivalent to $T$ for all $\varphi \in \text{M\"{o}b}$. 
\end{defn}
The problem of determining all the homogeneous normal operators, homogeneous contractions and homogeneous operators in the Cowen-Douglas class has been addressed in a series of papers \cite{BMsurvey,BHM, KM} previously.  
One of the goals of this paper is to determine modulo unitary equivalence, all hyponormal operators $T$ such that $[T^*,T]$ is in trace class that are homogeneous. This involves, among other things, finding a transformation rule for the Principal function of an operator under the M\"{o}bius transformations. 
\section{Preliminaries} 
An operator $T\in \mathcal L(\mathcal H)$ is said to be \textit{hyponormal} if the self-commutator $[T^*,T] = T^*T-TT^*$ is non-negative definite. 
A hyponormal operator $T$ is said to be a \textit{pure} if there is no nontrivial reducing subspace for $T$ on which it is normal. Every hyponormal operator $T$, modulo unitary equivalence, is of the form $T_p \oplus T_n$, 
where $T_p = T|_{\mathcal H_p}$, $T_n = T|_{\mathcal H_n}$ and 
$\mathcal H = \mathcal H_p \oplus \mathcal H_n$ such that $T_p$ is pure and $T_n$ is normal, see \cite[Theorem 1.3, Chapter II]{MP}. 

Any operator $T\in \mathcal L(\mathcal H)$ can be written in the form $T=A +iB$, 
where $A= \tfrac{T+ T^*}{2}$ and $B:= \tfrac{T- T^*}{2i}$ are self-adjoint. It follows that $[T^*,T] = 2 i[A,B]$. 

\subsection{Principal function} The Principal function of an operator $T$ is defined by means of an auxiliary operator valued function $E$ of two complex variables, called the determining function of $T$. The principal function $g_T$ of $T$ then appears by expressing the multiplicative determinant of the self-commutator, or the trace of the self-commutator $D:= [T^*,T]$ as an integral. We recall that the determining function $E$ is given by the formula 
\[
 E(z,w)= I-2 i D^{\tfrac{1}{2}} (A-z)^{-1} (B-w)^{-1} D^{\tfrac{1}{2}},\,\, z,w \in \mathbb C\setminus \sigma(A) \times \sigma(B).\]
Pincus in \cite{Pinc1,Pinc2} proved the existence of a function $g(v, u) \geqslant 0$ such that
\begin{equation}\label{detG}
 \det E(z, w)=\exp \left(\frac{1}{2 \pi i} \iint_\mathbb C g(u, v) \frac{d u}{u-z} \frac{d v}{v-w}\right),\,\, z,w \in \mathbb C\setminus \sigma(A) \times \sigma(B).
\end{equation}
The support of any almost everywhere determined version of $g(u,v)$ is said to be the ``determining set'' of the pair $A, B$, or equivalently that of the operator $T$. The essential closure of the determining
set is denoted by $D(A, B)$. It is proved in \cite{Pinc0} that $\sigma(T) = D(A,B)$. Thus, $\text{Supp}(g) \subseteq \sigma(T)$ and if $T$ is pure, then $\text{Supp}(g) = \sigma(T)$, see also \cite[$5^0$, pp. 105]{KC}.
\begin{rem}
 For every integrable, compactly supported function $g$ on $\mathbb C$, with $0 \leqslant g \leqslant 1$, there exists a pure semi-normal operator $T$, with $\operatorname{rank}\left[T^*, \mathrm{~T}\right]=1$ such that $[g]=\left[g_{T}\right]$ in $\mathrm{L}^{{1}}(\mathrm{~d} \mu)$.
The proof is in \cite[Theorem 1]{CPIUMJ}, see also \cite{Put}.
\end{rem}

\subsection{The tracial bi-linear form} 
Let $\mathbb C[x,y]$ denote the algebra polynomials over the complex field in the two indeterminates $x,y$.  Thus, any $p\in \mathbb C[x,y]$ is of the form 
\[p(x,y) = \sum_{j,k=1}^m a_{i,j} x^j y^k,\,\, a_{i,j}\in \mathbb C.\] 
Let $A, B$ be a pair of self adjoint operators in $\mathcal L(\mathcal H)$ such that $\|[A,B]\|_1 < \infty$.  Also, let $\mathbb C[A,B]$ be the algebra of operators generated by substituting 
$A, B$ in place of the commuting variables $x,y$ of the polynomial $p\in \mathbb C[x,y]$. Thus, if  $X,Y$ is any pair of operators in  $\mathbb C[A,B]$, 
 then the operator \[p(X,Y) = \sum_{j,k=1}^m a_{i,j} X^j Y^k\] 
is well defined modulo operators of trace class. The tracial bi-linear form associated with the pair $X,Y$ is  
\[(p,q) = \operatorname{tr} \,i\,[p(X,Y), q(X,Y)], p,q \in \mathbb C[x,y].\]
An amazing formula discovered by Helton and Howe \cite{HH} for the tracial bilinear form is given below. 
\begin{theorem*}[Helton-Howe] Suppose that $X,Y$ are a pair of operators such that $[X,Y]$ is in trace class. Then there exists a regular signed Borel measure
$\mu$ with compact support in $\mathbb C$ such that for $p,q\in \mathbb C[x,y]$, 
\[(p,q) = \operatorname{tr} \,i\,[p(X,Y), q(X,Y)] =  {\iint}_{\mathbb C} J(p,q) d\mu, \]
where $J(p,q)= \frac{\partial p}{\partial x} \frac{\partial q}{\partial y} - \frac{\partial p}{\partial y} \frac{\partial q}{\partial x}$. 
\end{theorem*}
Soon after the discovery of the Helton-Howe formula, Pincus established that the measure $\mu$ in the Helton-Howe formula is mutually absolutely continuous with respect to the area measure $dxdy$, that is,
$d\mu = g_T(x,y) dxdy$, where $g_T$ is the Principal function of the operator $T=A+iB$. 
 
\subsection{Unitary invariants} For $z,w$ in a neighbourhood of infinity, the operator valued determining function $E(z,w)$  of an irreducible pure hyponormal operator $T$ of trace class is a complete unitary invariant of $T$. The principal function $g_T$, on the other hand, is a unitary invariant in general but it is a complete invariant when the rank of $[T^*,T]$ is $1$. In what follows, we assume that the operator $T$ is a irreducible hyponormal (hence, pure), and that rank of $[T^*,T]=1$.  Thus, we assume without loss of generality that $[T^*,T] = x\otimes x$ for some $x\in \mathcal H$, where $x \otimes x$ denotes the \textit{non-negative definite} rank one operator $h \mapsto \langle h,x \rangle x$, $h\in \mathcal H$. In this case the multiplicative commutator and therefore, the determining function $E$ of the operator $T$ can be calculated explicitly as follows:
For any pair of complex numbers $z,w$ not in the spectrum of $T$, the operators $(T^*-\bar{w})^{-1}$ and $(T-z)^{-1}$ exist and the multiplicative commutator 
\begin{align*}
    (T-z)(T^*-\bar{w})(T-z)^{-1}(T^*-\bar{w})^{-1},
\end{align*}
is in the \textit{determinant class}, that is, it is of the form $I+K$, where $K$ is  trace class: 
\begin{multline*}
     (T-z)(T^*-\bar{w})(T-z)^{-1}(T^*-\bar{w})^{-1}\\ 
     = ((T^*T-x\otimes x-zT^*-\bar{w} T+z\bar{w})(T-z)^{-1}(T^*-\bar{w})^{-1})\phantom{FAN2FAN21}\\
     = ((T^*-\bar{w})(T-z)(T-z)^{-1}(T^*-\bar{w})^{-1} - (x\otimes x)(T-z)^{-1}(T^*-\bar{w})^{-1})\\
     = (I-(x\otimes x)(T-z)^{-1}(T^*-\bar{w})^{-1})\phantom{phantomphantomphantomphantom}\\
     =I + K,\phantom{phantomphantomphantomphantomphantomPphantomPphantomPphan} 
\end{multline*}
where $K= -(x\otimes x)(T-z)^{-1}(T^*-\bar{w})^{-1}$ is in trace class, and  $\operatorname{tr} K =  -\langle (T^*-\bar{w})^{-1}x, (T^*-\bar{z})^{-1}x \rangle$. Therefore, 
\begin{align*}
    \det(I-K) &= \exp\big (\operatorname{tr}\log(I-K)\big ) \\
    &= \exp\Big (\operatorname{tr}\sum_{j=1}^{\infty} \frac{(-1)^{j-1}}{j} (-K)^j\Big) \\
    &=\exp\Big (\sum_{j=1}^{\infty}\frac{(-1)^{2j+1}}{j}\operatorname{tr}(K^j)\Big)\\ &= \exp\Big (\sum_{j=1}^{\infty}\frac{(-1)^{2j+1}}{j}(\operatorname{tr}(K))^j\Big )\\ &=\exp\Big(\log(1-\operatorname{tr}(K))\Big)\\ &= 1-\operatorname{tr} (K) \\ &= 1 - \langle (T^*-\bar{w})^{-1}x, (T^*-\bar{z})^{-1}x \rangle 
\end{align*}
Therefore, combining with the formula \eqref{detG}, we have the equality:
\begin{align}\label{PincG}
    1 - \langle (T^*-\bar{w})^{-1}x, (T^*-\bar{z})^{-1}x \rangle  = \exp\left(-\frac{1}{\pi}\int_{\mathbb C}\frac{g_T(\zeta)}{(\zeta-z)(\bar{\zeta}-\bar{w})} dA(\zeta)\right).
\end{align}
For a different approach to establishing this formula, see \cite[Theorem 4.3]{CS}. 
We reiterate that the principal function of a irreducible hyponormal operators $T$ with rank-one self-commutator $x\otimes x$ is a complete unitary invariant of $T$. This remarkable theorem is due to Pincus and is in \cite{Pinc1}. A different unitary invariant is in \cite{CW}.


\subsection{An example}\label{shiftex} Let $S$ be the unilateral shift operator acting on the Hilbert space $\ell_2$ of square summable complex sequences by the rule: $S e_k = e_{k+1}$, where $\{e_0,e_1,e_2,...\}$ is the standard basis of $\ell_2$.
The self-commutator $[S^*,S] = e_0\otimes e_0$.  Since  $(S^*-\bar{w} I)^{-1}e_0 =- \frac{1}{\bar{w}}e_0$, we have that 
\begin{align*}
    1 - \langle (S^*-\bar{w} I)^{-1}e_0, (S^* - \bar{z} I)^{-1}e_0\rangle   = 1 - \langle -\frac{1}{\bar{w}}e_0, -\frac{1}{\bar{z}}e_0\rangle 
    & = 1 - \frac{1}{z\bar{w}}.
\end{align*}
We claim that 
\begin{align*}
    \exp\left(-\frac{1}{\pi}\int_{\overbar{\mathbb D}}\frac{1}{(\zeta - z)(\zb - \wb)}dA(\zeta)\right) = 1 - \frac{1}{z\bar{w}}.
\end{align*}
Taking  $|\zeta| \leq 1$, and $|z|,|w|>1$, and expanding $\frac{1}{\zeta - z}$ as well as $\frac{1}{\bar{\zeta} - \bar{w}}$ in a power series of $\frac{\zeta}{z}$ and $\frac{\bar{\zeta}}{\bar{w}}$, respectively, the claim is verified by integrating the product term by term. Thus the principal function of the unilateral shift $S$ is the characteristic function $\mathbb{1}_{\overbar{\mathbb D}}$ of the closed unit disc $\overbar{\mathbb D}$. 
\begin{rem} Let $\sigma_{\text {ess }}(T)$ be the essential spectrum 
of an operator $T$. 
For $\lambda \in \mathbb{C} \backslash \sigma_{\text {ess }} (T)$, the principal function $g(\lambda)=- \operatorname{ind} (T-\lambda)$, see \cite[$5^0$, pp. 105]{KC}.  Consequently, the principal function $g_S$ of the unilateral shift $S$ is $\mathbb 1_{\overbar{\mathbb D}}$. 
\end{rem}
\section{The action of the M\"{o}bius group} 
The hyponormal operators share an important property with normal operators, namely, the spectral radius $\rho(T)$ of a hyponormal operator equals its norm $\|T\|$. However, unlike normal operators, if $T$ is a \textit{pure} hyponormal operators, then by Putnam's inequality, the area measure of spectrum $\sigma(T)$ must be 
positive.  
\subsection{Invariance} It is not hard to verify that if $T$ is hyponormal, then $\varphi(T)$ is also hyponormal for any $\varphi$ in M\"{o}b, the biholomorphic automorphism group of the unit disc $\mathbb D$. We reproduce the proof below from \cite[Lemma 1]{JS}. 
\begin{proposition}[Stampfli]\label{Stamp}
If $T$ is hyponormal, then $\varphi(T)$, $\varphi$ in M\"{o}b, is also hyponormal. 
\end{proposition}
\begin{proof}
Any M\"{o}bius transformation is a composition of an affine transformation and an inversion of some other affine transformation. We have 
\begin{align*}
    [(aT + b)^*,aT + b)] = |a|^2[T^*,T] \geq 0.
\end{align*}
Therefore, to complete the proof, it is enough to verify that $[(T^*)^{-1}, T^{-1}]$
is hyponormal. By hypothesis, we have that   
\begin{align*}
    0 \leqslant T^{-1}(T^*T-TT^*)(T^*)^{-1} = T^{-1}T^*T(T^*)^{-1} - I
\end{align*}
If $A$ is invertible and $A \geqslant I$, then $A^{-1} \leqslant I$.  Therefore,
\begin{align*}
    I - T^*T^{-1}(T^*)^{-1}T \geqslant 0.
\end{align*}
Hence, \[[(T^*)^{-1},T^{-1}] = ((T^*)^{-1}T^{-1} - T^{-1}(T^*)^{-1}) = (T^*)^{-1}(I - T^*T^{-1}(T^*)^{-1}T)T^{-1} \geqslant 0\]
completing the proof of the proposition.
\end{proof}
We now re-write the formula for the tracial bi-linear form in complex co-ordinates and in slightly greater generality, see \cite[Chapter X, Theorem 2.4, and  Equation (12), pp. 242]{MP}. 
\begin{thm}[Carey-Helton-Howe-Pincus] Suppose that $T\in \mathcal L(\mathcal H)$ is a hyponormal operator with $[T^*,T]$ is in the trace class $\mathcal S_1(\mathcal H)$. Then for any pair of functions $p,q$ in the Frechet Space $C^\infty(\sigma(T))$ of all smooth functions on $\sigma(T)$, we have the equality
\begin{align*}
    \operatorname{tr}[p(T,T^*), q(T,T^*)] = \frac{1}{\pi}\int_{\sigma(T)}J(p,q)g_T d\mu,
\end{align*}
where $J(p,q): = \tfrac{\partial{p}}{\partial\overbar{z}} \tfrac{\partial q}{\partial z}- \tfrac{\partial p}{\partial z} \tfrac{\partial q}{\partial \overbar{z}}$. 
\end{thm}
The proof of the following lemma follows directly from the Carey-Helton-Howe-Pincus formula.  

\begin{lemma} \label{lem:3.1} Suppose that $T\in \mathcal L(\mathcal H)$ is a hyponormal operator and  $[T^*,T]$ is in the trace class $\mathcal S_1(\mathcal H)$ and that $\sigma(T)\subseteq \overbar{\mathbb D}$.  Then  $[(T^*-\overbar{\lambda})^{-1},(T-\lambda)^{-1}]$ is also in $\mathcal S_1(\mathcal H)$ for $\lambda\notin\sigma(T)$. In particular, $[\varphi(T)^*,\varphi(T)]$ is in $\mathcal S_1(\mathcal H)$ for any $\varphi\in \text{M\"{o}b}$.
\end{lemma}

\begin{proof} Pick $p(\zeta,\zb) = \frac{1}{\zb-\lc}$ and $q(\zeta,\zb) = \frac{1}{\zeta-\lambda}$. Then $J(p,q) = \frac{1}{|\zeta-\lambda|^4}<k$ for some $k>0$ since $\lambda\notin\overbar{\mathbb D}$. 
Therefore,
\begin{align*}
    \ttr[(T^*-\overline{\lambda})^{-1},(T-\lambda)^{-1}] & = \ttr[p(T,T^*),q(T,T^*)]\\
    & = \frac{1}{\pi}\int_{\sigma(T)}\frac{1}{|\zeta-\lambda|^4}g_T  d\mu\\
    & \leqslant \frac{k}{\pi}||g_T||_{L^1(\sigma(T))} \\
    & < \infty. 
\end{align*}
Since affine transform of a trace class operator is again in trace class, the proof is complete. 
\end{proof}
We now compute the self commutator of the operator $\varphi(T)$. For this, we note that $\varphi(z) = \frac{z-a}{1-\bar{a}z} = -(\bar{a})^{-1}+c(z-\bar{a}^{-1})^{-1}$, where $c = \frac{a-\bar{a}^{-1}}{\bar{a}}$
 \begin{align*}
     [\varphi(T)^*,\varphi(T)] & = [-a^{-1} I+\bar{c}(T^*-a^{-1})^{-1}, -\bar{a}^{-1} I+c(T-\bar{a}^{-1})^{-1}]\\
     & = (-a^{-1} I+\cb(T^*-a^{-1})^{-1})(-\ab^{-1} I+c(T-\ab^{-1})^{-1})  \\&\phantom{SagarGhosh} - (-\ab^{-1} I+c(T-\ab^{-1})^{-1})(-a^{-1} I+\cb(T^*-a^{-1})^{-1})\\
     & = |c|^2((T^*-a^{-1})^{-1}(T-\ab^{-1})^{-1} - (T-\ab^{-1})^{-1}(T^*-a^{-1})^{-1})\\
     & = |c|^2 \big ((T-\ab^{-1})(T^*-a^{-1})\big )^{-1}  [(T^*-a^{-1}), (T-\ab^{-1})]\big ( (T^*-a^{-1}) (T-\ab^{-1}) \big )^{-1}\\ 
     & = |c|^2 \big ((T-\ab^{-1})(T^*-a^{-1})\big )^{-1}  [T^*,T] \big ( (T^*-a^{-1}) (T-\ab^{-1}) \big )^{-1}.
 \end{align*}
The computation of $[\varphi(T)^*,\varphi(T)]$ facilitates the proof of the lemma below. 
\begin{lemma}\label{lem:3.2}
     Suppose that $T\in \mathcal L(\mathcal H)$ is a hyponormal operator and  the rank of $[T^*,T]$ is $1$ and that $\sigma(T)\subseteq \overbar{\mathbb D}$. Then the rank of the self-commutator $[\varphi(T)^*, \varphi(T)]$ is also $1$. 
\end{lemma}
\begin{proof} For the proof, in view of the preceding discussion, it is enough to verify that whenever $T$ is an invertible operator with $\operatorname{rank} [T^*,T]=1$, the rank of $[{T^*}^{-1},T^{-1}]$ is also $1$. 
     By hypothesis, $[T^*,T] = x \otimes x$ for some vector $x$ in $\mathcal H$. Hence,
\begin{align*}
    T^{-1}T^*T(T^*)^{-1} & = I + T^{-1}(T^* T - TT^*)(T^*)^{-1}\\
     & = I + T^{-1}(x\otimes x)(T^*)^{-1}.
     \end{align*} 
 Taking inverses on both sides, we have 
 \begin{align*} 
     T^{-1}(T^*)^{-1} & = (T^*)^{-1}[I+T^{-1}(x\otimes x)(T^*)^{-1}]^{-1}T^{-1}\\
     & = (T\{I+ T^{-1}(x\otimes x)(T^*)^{-1}\}T^*)^{-1}\\
     & = (TT^*+x\otimes x)^{-1}.
\end{align*}
Similarly, 
\[(T^*)^{-1}T^{-1} = (T^*T - x\otimes x)^{-1}.\]
Therefore,
\begin{align*}
    [(T^*)^{-1}, T^{-1}] & = (T^*)^{-1}T^{-1} - T^{-1}(T^*)^{-1}\\
    & = (T^*T - x\otimes x)^{-1} - (TT^*+x \otimes x)^{-1}\\
    & = (T^*T - x\otimes x)^{-1}(\{TT^*+x\otimes x\} - \{T^*T-x\otimes x\})(TT^*+x\otimes x)^{-1}\\
    & = [T^*T - x\otimes x]^{-1}[x\otimes x][TT^*+ x\otimes x]^{-1}\\
    & = (TT^*)^{-1}(x\otimes x)(T^*T)^{-1}.
\end{align*}
It follows that the self-commutator of $T^{-1}$ is also of rank one completing the proof. 
\end{proof}
\begin{rem} Combining Proposition \ref{Stamp} and Lemma  \ref{lem:3.2}, we conclude that the set of \textit{pure} hypnormal operators with rank $1$ self-commutator is left invariant under the action of the M\"{o}bius group. Similarly, combining Proposition \ref{Stamp}, this time with Lemma \ref{lem:3.1}, we see that the set of \textit{pure} hypnormal operators $T$ with $\|[T^*,T\|_1$ finite is also left invariant under the action of M\"{o}b.
\end{rem}
\subsection{A change of variable formula for the principal function}
A change of variable formula for the principal function appears in \cite[pp. 106 - 107]{KC} and also in \cite[pp. 245]{MP}. However, for our purposes, we need a change of variable formula for the principal function in the form given below.
\begin{proposition} 
Let $T$ be a pure hyponormal operator with trace class self-commutator and set $W:=\varphi(T)$, $\varphi$ in M\"{o}b. Assume that the spectrum of $T$ is contained in the closed unit disc.   Then the relationship between the two principal functions $g_T$ and $g_W$ is given by the change of variable formula 
\begin{align*}
    g_W(\zeta) = g_T(\varphi^{-1}(\zeta)), \zeta \in \sigma(W).
\end{align*}
\end{proposition}

\begin{proof} We have proved that $W$ is a hyponormal operator with $\|[W^*,W]\|_1 < \infty$. We note that  $\varphi(T)^* = \varphi^*(T^*)$, where $\varphi^*(z) = \overline{\varphi(\bar{z})}$. Setting $\tilde{p}(z,\bar{z}):=p(\varphi(z),\overline{\varphi(z)})$ and  $\tilde{q}(z,\bar{z}):=q(\varphi(z),\overline{\varphi(z)})$, we have that 
\begin{align*}
    \operatorname{tr}[p(\varphi(T),\varphi(T)^*),q(\varphi(T),\varphi(T)^*)] = \frac{1}{\pi}\int_{\sigma(\varphi(T))}J(p,q)g_{\varphi(T)}(\zeta) dA(\zeta).
\end{align*}
On the other hand, 
\begin{align*}
    \operatorname{tr}[p(\varphi(T),\varphi(T)^*),q(\varphi(T),\varphi(T)^*)] & = \operatorname{tr}[\pp(T,T^*),\qq(T,T^*)]\\
    & = \frac{1}{\pi}\int_{\sigma(T)}J_\zeta(\pp,\qq)g_T(\zeta)dA(\zeta)\\
    & = \frac{1}{\pi}\int_{\sigma(\varphi(T))}J_\eta(p,q)g_T(\varphi^{-1}(\eta))dA(\eta),
\end{align*}
where $\eta = \varphi(\zeta)$. By the chain rule, we have 
$\frac{\partial \tilde{p}}{\partial \bar{\zeta}} = \frac{\partial \tilde{p}}{\partial \bar{\eta}} \frac{\partial \varphi}{\partial \bar{\zeta}}$, and  similarly $\frac{\partial \tilde{p}}{\partial {\zeta}} = \frac{\partial \tilde{p}}{\partial {\eta}} \frac{\partial \varphi}{\partial {\zeta}}$.
Thus, we have the equality  
\begin{align*}
    J_{\zeta}(\tilde{p},\tilde{q}) = J_{\eta}(p,q)\Big (\frac{\partial(\overline{\varphi(\zeta))}}{\partial{\bar{z}}}\frac{\partial(\varphi(\zeta))}{\partial{\zeta}}\Big ).
\end{align*}
Consequently, 
\begin{align*}
    dA(\eta) & = -\frac{1}{2i}d\eta\wedge d\overline{\eta}\\
    & = -\frac{1}{2i}\Big (\frac{\partial(\overline{\varphi(\zeta))}}{\partial{\zb}}\frac{\partial(\varphi(\zeta))}{\partial{\zeta}}\Big)d\zeta\wedge d\zb\\
    & = \Big (\frac{\partial(\overline{\varphi(\zeta))}}{\partial{\zb}}\frac{\partial(\varphi(\zeta))}{\partial{\zeta}}\Big )dA(\zeta).
\end{align*} 
Hence,
\begin{align*}
    J_\zeta(\pp,\qq) dA(\zeta) = J_\eta(p,q)dA(\eta).
\end{align*}
Since $p$ and $q$ are arbitrary $C^\infty$ functions on $\sigma(T)$, we conclude that 
\begin{align*}
    g_{\varphi(T)}(\zeta) = g_T(\varphi^{-1}(\zeta))
\end{align*}
completing the proof. 
\end{proof}
\section{Homogeneous hyponormal operators $T$ with $\operatorname{rank} [T^*,T] =1$}
We have already remarked that the principal function of a pure hyponornmal operator in the trace class $\mathcal S_1(\mathcal H)$ is not a complete unitary invariant for the operator $T$ in general. However, it is not hard to see that it is a unitary invariant.
\begin{proposition} Let $T$ be a pure hyponormal operator in $\mathcal S_1(\mathcal H)$. If $W$ is an operator unitarily equivalent to $T$, then the principal functions of $W$ and $T$ coincide. 
\end{proposition} 
\begin{proof} Let $W = UTU^*$ for some unitary operator $U$. conjugation of $T$. The operator $W$ is hyponormal and is in $\mathcal S_1(\mathcal H)$.
For any polynomial $p\in \mathbb C[z,\bar{z}]$, we have $p(W,W^*) = Up(T,T^*)U^*$. Hence, by the Helton-Howe formula, we find that 
\begin{align*}
    \frac{1}{\pi}\int_{\sigma(W)}J(p,q)g_W(\zeta)dA(\zeta) & = \operatorname{tr} [p(W,W^*),q(W,W^*)]\\
    & = \operatorname{tr}  [Up(T,T^*)U^*,Uq(T,T^*)U^*]\\
    & = \operatorname{tr}  (U[p(T,T^*),q(T,T^*)]U^*)\\
    & = \operatorname{tr} [p(T,T^*),q(T,T^*)]  \\
    & = \frac{1}{\pi}\int_{\sigma(T)}J(p,q)g_T(\zeta)dA(\zeta)
\end{align*} 
Since $\sigma(T)= \sigma(W)$, we have that 
$\frac{1}{\pi}\int_{\sigma(T)}J(p,q)(g_T - g_W)(\zeta)dA(\zeta) = 0$, $p,q \in \mathbb C[x,y]$, and in consequence $g_W = g_T$. 
\end{proof} 
Imposing the condition of homogeneity on a pure hyponormal opeartor $T\in \mathcal S_1(\mathcal H)$, we investigate what happens to the principal function $g_T$. 

We begin with the simple observation that if $T$ is a homogeneous operator, then by the spectral mapping theorem, the spectrum $\sigma(T)$ must be invariant under the action of the M\"{o}bius group. Consequently, $\sigma(T)$ has to be either the closed unit disc $\overbar{\mathbb D}$, or the unit circle $\mathbb T$. However, if $T$ is also a pure hyponormal operator, then as we have noted earlier, $\sigma(T)$ cannot be $\mathbb T$. What is more, 
\begin{proposition} \label{prop:4.2}
Suppose that $T$ is a pure hyponormal homogeneous operator such that $[T^*,T]$ is in $\mathcal S_1(\mathcal H)$. Then the principal function $g_T$ is constant on the spectrum $\sigma(T)$.  
\end{proposition} 
\begin{proof} Since $\varphi(T)$ is unitarily equivalent to $T$, $\varphi$ in M\"{o}b, it follows that $ g_T(z) = g_{\phi(T)}(z)$. By the change of variable formula for the principal function, we have 
$g_{\phi(T)}(z) = g_T(\phi^{-1}(z))$. Combining these two equalities, we conclude that 
\begin{equation}\label{eqn:4.1}
    g_T(z) = g_{\varphi(T)}(z) = g_T(\varphi^{-1}(z)), 
\end{equation}
for all $\varphi \in \mbox{M\"{o}b}$. For a fixed but arbitrary $z\in \mathbb D$, pick a M\"{o}bius transformation $\varphi_z$ with the property: $\varphi_z(0)=z$. 
Using this $\varphi_z$ in Equation \eqref{eqn:4.1}, we have
\begin{align*}
    g_T(z) = g_{\phi_z(T)}(z) = g_T({\phi_z}^{-1}(z)) = g_T(0).
\end{align*}
We therefore conclude that $g_T$ must be a constant on $\sigma(T)$ with $0< g_T(0) \leq 1$.  
\end{proof} 
We have now all the tools to prove the only new result of this short note. 
\begin{thm}
  The only homogeneous \textit{pure} hyponormal operator $T$ with $\operatorname{rank}[T^*,T] = 1$, modulo unitary equivalence, is the unilateral shift.  
\end{thm}
\begin{proof}
We have already shown that the principal function $g_S$ of the unilateral shift $S$ is constant on the spectrum $\overbar{\mathbb D}$ of $S$. Indeed, $g_S(z) = 1$, $z\in \overbar{\mathbb D}$.

We have also shown that the principal function $g_T$ of a homogeneous pure hyponormal operator $T$ with $\operatorname{rank}[T^*,T] = 1$ must be a constant and moreover, $0 < g_T \leqslant 1$. So, to complete the proof, we have to simply show that there is no such operator with $g_T = c< 1$. Let us suppose to the contrary that there exists such an operator $T$ with $g_T = c < 1$. 
In the determinant expansion formula  \eqref{PincG}, setting $g_T=c$, we have (as in Example \ref{shiftex}): 
\begin{align}\label{align:0}
    1-\langle (T^*-\bar{w})^{-1}x,(T^*-\bar{z})^{-1}x \rangle &= \exp\Big (-\frac{1}{\pi}\int_{\sigma(T)}\frac{g_T(\zeta)dA(\zeta)}{(\zeta-z)(\bar{z}-\bar{w})}\Big )\\ &= \Big (1-\frac{1}{z\bar{w}}\Big )^c.
\end{align}
Putting $z=w$ in Equation \eqref{align:0} 
we have the equality 
\begin{equation} \label{eqn:4.2}  1 - \|(T^*-\wb)^{-1}x\|^2 = \Big(1-\frac{1}{|w|^2}\Big)^c. 
\end{equation}
Since $T$ is homogeneous and hyponormal, the spectrum $\sigma(T)$ can only be $\overbar{\mathbb D}$, the possibility of $\sigma(T) = \mathbb T$ is ruled out by Putnam's inequality. 
For a hyponormal operator, the spectral radius $\rho(T)  = \|T\|$ and we conclude that   that $\|(T-wI)^{-1}\| = \rho{((T-wI)^{-1})} \leqslant \frac{1}{|w|}$. Since $[T^*,T] = x\otimes x$ for some $x\in \mathcal H$ by hypothesis, taking $p(z,\bar{z}) = \bar{z}$ and $q(z,\bar{z}) = z$ in the Helton-Howe formula we conclude that $\|x\| = \sqrt{c}$.  Hence, the inequality $\|(T^*-\wb)^{-1}x\|\leqslant \sqrt{c}\|(T^*-\wb)^{-1}\|$, gives
\begin{align} \label{align:4.3}
    1 - \|(T^*-\wb)^{-1}x\|^2 \geq 1 - \|(T^*-\wb)^{-1}\|^2\|x\|^2 
    = 1 - c\|(T^*-\wb)^{-1}\|^2 
    \geqslant  1 - \frac{c}{|w|^2}.
\end{align}
Combining the equality \eqref{eqn:4.2} with the inequality \eqref{align:4.3}, we have 
\begin{align} \label{align:4.4}
    \Big(1-\frac{c}{|w|^2}\Big ) \leqslant \Big (1-\frac{1}{|w|^2}\Big )^c, \,\, |w|>1.
\end{align} 
It is easy to verify that the inequality \eqref{align:4.4} is false unless $c=1$ completing the proof. 
\end{proof}
We now give example of a large class of unitarily inequivalent operators possessing the same principal function. Also, see reamrk below Lemma 1 in \cite[pp. 252]{putex}.

Let $\{T_{\lambda}\}_{\lambda>1}$ be the weighted shift operator  with weight sequences $\{w_n(\lambda)\}_{n\geqslant 0}$,  $w_n(\lambda) = \tfrac{n+1}{n+\lambda}$. For $\lambda>1$, the weight sequence $\{w_n(\lambda)\}$ is strictly increasing and hence $T_\lambda$ is hyponormal. The operator $T_\lambda$ is also pure and cyclic. 
Moreover, 
\begin{align*}
    \operatorname{tr}[T_\lambda^*,T_\lambda] & = \sum_{i=0}^{\infty}(w^2_{i+1}(\lambda)-w^2_{i}(\lambda))+w^2_0(\lambda) =1.
\end{align*}
For $\lambda_1\neq\lambda_2$, the two operators $T_{\lambda_1}$ and $T_{\lambda_2}$ are unitarily inequivalent. But all these operators are homogeneous, see \cite{BMJFA}. Therefore, the principal function  $g_{T_\lambda}$ is constant, say $c$, on $\overbar{\mathbb D}$. But then
\begin{align*}
  1= \operatorname{tr}[T_\lambda^*,T_\lambda] = \frac{1}{\pi}\int_{\overbar{\mathbb D}} c\,dA(\zeta). 
\end{align*}
Thus, $c=1$ and it follows that $g_{T_\lambda}$ is identically $1$ on $\overbar{\mathbb D}$ for all $\lambda>1$.
\subsection{Open problem} Find all the pure hyponormal operators $T$ such that $[T^*,T]$ is in $\mathcal S_1(\mathcal H)$ and that 
$g_T$ is constant on $\sigma(T)$ modulo unitary equivalence. 
\begin{rem}
    In studying homogeneous contractions $T$ assuming that both the defect indices of $T$ are equal to $1$, it was shown that the Sz.-Nagy--Foias characteristic function of $T$ must be constant. This observation leads to  a class of homogeneous bi-lateral shifts (all of them inequivalent among themselves), parametrized by $c >0$,  possessing a constant characteristic function, see \cite{CM, const}. 

Similarly, homogeneous operators $T$ in the Cowen-Douglas class $B_1(\mathbb D)$ are determined by specifying the curvature $\lambda = -\mathcal K_T(0) > 0$ just at one point. From this, one infers that an operator $T$ in $B_1(\mathbb D)$ is homogeneous if and only if $T$ is of the form $T_\lambda^*$, $\lambda >0$, discussed above (see \cite{GM}).

    The situation involving the hyponormal operators $T$ with $\operatorname{rank} [T^*,T] = 1$, appears to be very different. 
    Here again, the unitary invariant $g_T$, under the assumption of homogeneity, is a constant function, say $c$, with $0 < c \leqslant 1$. But there is only one \textit{homogeneous} hyponormal operator $T$ with $[T^*,T] =x\otimes x$, namely, the unilateral shift corresponding to $c=1$.  
\end{rem}
\subsubsection*{Acknowledgment} We are very grateful to K. B. Sinha for his generous help during the preparation of this manuscript. In particular, the proof of Lemma \ref{lem:3.2} evolved during discussions with him. We also express our gratitude to M. Putinar for patiently answering some of our questions and pointing to the paper \cite{putex}. 
\subsubsection*{Postscript} In a conversation with the second author, in the year 1983, Kevin Clancey had remarked that the only homogeneous \textit{pure} hyponormal operator with rank $1$ self-commutator might be the unilateral shift. We have verified this statement to be correct in this short note.

\end{document}